# A Bilateral Reserve Market for Variable Generation: Concept and Implementation

Runze Chen, *Student Member, IEEE*, Audun Botterud, *Member, IEEE*, Hongbin Sun, *Senior Member, IEEE*, Yang Wang, *Student Member, IEEE*

*Abstract*—Substantial changes in the generation portfolio take place due to the fast growth of renewable energy generation, of which the major types such as wind and solar power have significant forecast uncertainty. Reducing the impacts of uncertainty requires the cooperation of system participants, which are supported by proper market rules and incentives. In this paper, we propose a bilateral reserve market for variable generation (VG) producers and capacity resource providers. In this market, VG producers purchase bilateral reserve services (BRSs) to reduce potential imbalance penalties, and BRS providers earn profits on their available capacity for re-dispatch. We show in this paper that by introducing this product, the VG producers' overall imbalance costs are linked to both their forecast quality and the available system capacity, which follows the cost-causation principle. Case studies demonstrate how the proposed BRS mechanism works and its effectiveness.

## NOMENCLATURE

### A. Indices, parameters, sets and functions

| | |
|---|---|
| $p^v$ | Actual output of a VG producer. |
| $P^{v,max}$ | Installed capacity of a VG producer. |
| $\alpha^+, \alpha^-$ | Penalty factors for over- and under-generation. |
| $\lambda^D$ | Day-ahead energy price. |
| $\lambda^R$ | Real-time energy price. |
| $\pi^+, \pi^-$ | Prices of upward and downward BRS. |
| $f(\cdot)$ | Probability density function of VG output. |
| $F(\cdot)$ | Cumulative probability density function of VG output. |

### B. Variables

| | |
|---|---|
| $\hat{p}^v$ | Day-ahead schedule of a VG producer. |
| $\Delta r$ | Executed amount of BRS. |
| $r^+, r^-$ | Contract amount of downward and upward BRS. |
| $R$ | Revenue of a VG producer without BRS. |
| $R^*$ | Revenue of a VG producer with BRS. |
| $R^g$ | Revenue of a dispatchable unit without BRS. |
| $R^{g*}$ | Revenue of a dispatchable unit with BRS. |

## I. INTRODUCTION

Due to economic and environmental concerns, the world is now experiencing a trend of deepening renewable energy penetration. Most common types of renewable energy, such as wind and solar power, are also named as variable generation (VG) because of their nature of variability and uncertainty. Proportionally, in some countries such as Denmark and Spain the share of wind and solar power in electricity production has reached over 30% [1], and the rate is still increasing rapidly. These significant changes pose many challenges to power system and market operation.

From the system perspective, sufficient amount of flexibility should be ensured to maintain operational reliability. Currently, there is consensus that the rising share of VG has increased the requirements of operating reserves [2]. Moreover, considering the nature of forecasting errors, it is found more economically efficient to allocate reserve needs dynamically, i.e. quantifying reserve needs for different time instants according to varying levels of uncertainty [3]-[5].

However, operating reserves are mainly kept for short-term variability and uncertainty. Capacities reserved for re-dispatch, which address the deviations between day-ahead (DA) and real-time (RT) forecasts, are not accurately defined and compensated by the current market rules. Considering that the increasing integration of VG will cause more deviations between DA schedule and RT dispatch, insufficient re-dispatch capacity might lead to more severe fluctuating RT prices, which is undesirable for market operation. Instead of incenting the provision of re-dispatch capacity only through varying RT prices, a promising way is to regard them as a reserve service, and compensate them independently from the energy market.

Regardless of the specific market rules, VG will bring additional operation costs to the system. Traditionally, these integration costs are allocated to consumers based upon energy consumption volumes. This seems to violate the rule of cost causation because the VG producers bear little responsibility for the additional costs they incur on the system. Actually, some markets, especially those with large share of VGs, have begun to put more responsibility to the VG producers. In Nord Pool countries, Spain, UK and Netherland, VG producers are now exposed to imbalance costs if they deviate from DA schedules [6]. In MISO, VG producers (DIRs) are also subject to additional charges for deviating from DA settlements [7]. However, an accurate allocation of integration cost is still difficult because a widely acceptable method to calculate the cost is hardly available [8]. Therefore, a market that addresses both the value of re-dispatch capacity and the responsibility of managing uncertainty might be a potential alternative. In the scheme proposed in this paper, the compensation of available capacity and allocation of costs are done in the market process without centralized decisions.

The basic idea of the market is that the provision of re-dispatch capacity is defined as a bilateral reserve service (BRS), and that the VG producers can purchase BRS to partially avoid potential deviation penalties. BRS is very similar to traditional reserve products, but has a longer time horizon and is carried out in a bilateral way. Actually, bilateral capacity contract is not an entirely new concept. In PJM, load serving entities (LSE) can fulfill their capacity requirements by purchasing it from other dispatchable resources [9]. However, this cannot be directly applied to VG producers because this is mainly for long-term capacity adequacy, and more importantly, the capacity obligation of a LSE is clearly defined while that of a VG producer is not. BRS for VG producers has already been discussed by [10] and [11], but the authors mainly focused on the pricing issues and did not explain the market design of such products. Further analysis is still needed to understand the potential benefits of this mechanism in current markets.[1]

In this paper, we will discuss the market implementation for BRS, and analyze the demand characteristics of this service. The remainder of this paper is organized as follows: in Section II, we describe the concepts of BRS and the market implementation; in Section III and IV, we analyze the demand characteristics; in Section V, the impacts of network constraints are briefly discussed. Case studies are presented in Section VI, followed by conclusions in Section VII.

## II. Bilateral Reserve Services For VG Producers

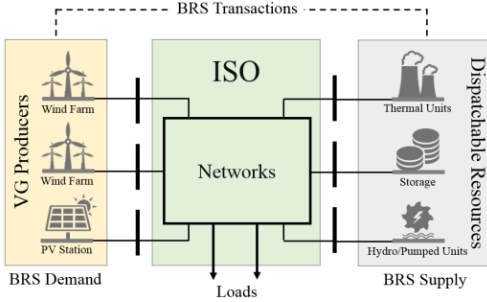

Fig. 1. Simplified structure of power system and market with BRS transactions.

Fig. 1 presents a simplified structure of power system and market. In terms of bilateral transactions, VG producers are on the demand side and dispatchable resources serve as providers of BRS. To establish the responsibility of managing uncertainty, VG producers are subject to imbalance penalties if they deviate from DA schedules. These penalties can be partially avoided by purchasing BRS.

Fig. 2 demonstrates how a BRS contract works. To own the BRS, the producer has to pay premiums equal to the sum of payments for upwards and downward BRS, i.e. $\pi^+ r^+ + \pi^- r^-$. Approaching RT, when the VG producer has higher certainty about its realized output level, it can choose to execute the BRS within purchased capacity to avoid at least a part of its imbalance costs. For example, in Scenario 1, the over-generation is within the purchased amount of downward BRS. The VG producer will therefore execute a part of its downward BRS to modify its DA schedule to eliminate or reduce the day-ahead and real-time deviation. Correspondingly, the BRS provider's DA schedule will be modified in the opposite direction at the same amount. For example, if the DA schedule of the VG producer is 100 MW and that of the BRS provider is 200 MW, and a 20 MW downward BRS is executed. Then, the DA schedule of the VG producer and BRS provider will be modified to 120 MW and 180 MW respectively, leaving the total amount, i.e. 300 MW, unchanged. Similarly, Scenario 2 describes the execution of upward BRS. In Scenario 3 & 4, the deviations exceed the BRS capacity. Then, the VG producer will choose to execute all BRS but only a portion of the deviations can be covered.

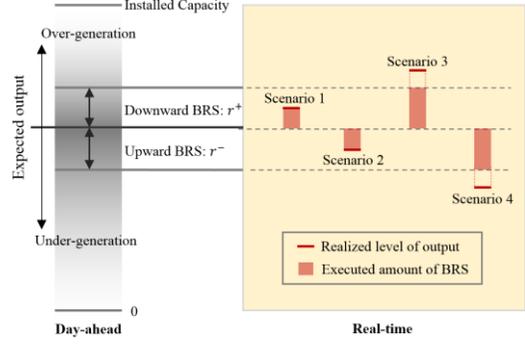

Fig. 2. The purchasing and execution of BRS for a VG producer.

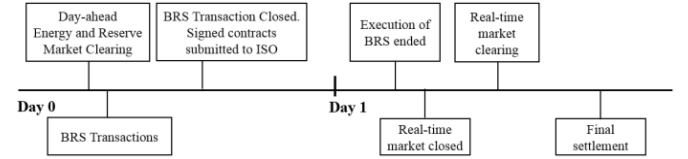

Fig. 3. BRS transaction and execution along with DA and RT market (the processes in Day 1 are repeated for each time instant).

The transaction and execution of BRS can be embedded into the DA and RT market structure. As presented in Fig. 3, the BRS transactions are open for a certain period of time after the DA energy and reserve market is cleared. The amount of BRS sold by a unit is constrained by its physical capacity, as well as scheduled energy and reserves. After the transaction period is closed, all the signed BRS contracts are sent to the ISO for validation. Then, before the closing of each RT market, the VG producers have to claim their execution decision, and the unused BRS will be released. Finally, the settlement of VG producers and BRS providers will be based on both DA and RT market settlements, as well as all the executed BRS contracts.

Executing a BRS contract is intrinsically switching the DA positions of producers, and may change the DA power flow. If the buyer and seller of BRS are located in the same node, the other players in the system will be unaffected. However, when they are at different nodes, the network constraints need to be carefully examined. In this paper, we ignore the network constraints, but include a brief discussion in Section IV addressing this issue.

BRS provides the connection between the sources of uncertainty and the providers of re-dispatch capacity. It encourages dispatchable resources to provide capacity-based services and thereby increase their profits. Moreover, the BRS market is embedded with the information about uncertainty and cost of capacity, based on which a market-based price for BRS is automatically decided. Next, we will explore the features of this market in detail.

## III. Demand and Supply Characteristics of BRS

### A. Demand

As mentioned previously, though not universally applicable, it is assumed that VG producers are subject to imbalance penalties if they deviate from DA schedules. Otherwise, the VG producers bear no responsibility for their forecast errors and there is no need for subsequent discussions. Current market practices vary with regards to the formation of deviation penalties. It might be a fixed dollar amount, or relate to the DA or RT clearing prices. Certainly, we are unable to explore all the possible penalty mechanism. But, proper assumptions can be made without loss of generality.

---

[1] In this paper we focus on BRS for VG, but other market participants (e.g. loads) could also benefit from purchasing the same product.

In this paper, it is assumed that VG producers have to sell its over-generation at lower prices and settle its under-generation at higher prices than what occurs in the DA market. The pay-off function of a VG producer can be expressed as:

$$R = \begin{cases} \lambda^D \hat{p}^v + (1-\alpha^+)\lambda^D (p^v - \hat{p}^v), & \hat{p}^v \leq p^v \leq P^{v,\max} \\ \lambda^D \hat{p}^v - (1+\alpha^-)\lambda^D (\hat{p}^v - p^v), & 0 \leq p^v \leq \hat{p}^v \end{cases} \quad (1)$$

The punitive RT settlement prices are represented by proportional relationships to the DA market clearing prices. In practice, if penalties are given and fixed at the DA stage the penalty factors $\alpha^+$ and $\alpha^-$ can explicitly calculated. If they are uncertain until the realization of RT clearing prices, the factors can be estimated based on historic data and forecast information [12].

If the producers are allowed to sign BRS contracts, the revenue will be as

$$R^* = \begin{cases} \lambda^D (\hat{p}^v + r^+)\hat{p}^v + (1-\alpha^+)\lambda^D (p^v - \hat{p}^v - r^+), & \hat{p}^v + r^+ \leq p^v \leq P^{v,\max} \\ \lambda^D p^v, & \hat{p}^v - r^- \leq p^v \leq \hat{p}^v + r^+ \\ \lambda^D (\hat{p}^v - r^-) - (1+\alpha^-)\lambda^D (\hat{p}^v - r^- - p^v), & 0 \leq p^v \leq \hat{p}^v - r^- \end{cases} \quad (2)$$

Therefore, the revenue becomes a function of purchased BRS $r^+$ and $r^-$. To determine the optimal amount of BRS, a VG producer aims at maximizing its expected revenue, which can be expressed as

$$\max \mathbb{E}[R] = \int_0^{\hat{p}^v - r^-} \left[ \lambda^D (\hat{p}^v - r^-) - \lambda^D (1+\alpha^-)(\hat{p}^v - r^- - p^v) \right] f(p^v) dp^v$$
$$+ \int_{\hat{p}^v - r^-}^{\hat{p}^v + r^+} \lambda^D p^v f(p^v) dp^v \quad (3)$$
$$+ \int_{\hat{p}^v + r^+}^{P^{v,\max}} \left[ \lambda^D (\hat{p}^v + r^+) + \lambda^D (1-\alpha^+)(p^v - \hat{p}^v - r^+) \right] f(p^v) dp^v$$

By taking the partial derivative of $r^+$ and $r^-$, we have:

$$\frac{\partial \mathbb{E}[R]}{\partial r^+} = \lambda^D \hat{\alpha}^+ \left[ 1 - F(\hat{p}^v + r^+) \right] \quad (4)$$

$$\frac{\partial \mathbb{E}[R]}{\partial r^-} = \lambda^D \hat{\alpha}^- F(\hat{p}^v - r^-) \quad (5)$$

Equation (4) and (5) describe the decreasing marginal utility of downward and upward BRS, i.e. the demand characteristics of BRS. According to the optimality condition, given the price level of downward and upward BRS, the optimal amount of BRS can be determined as:

$$r_0^+ = F^{-1}(1 - \pi_0^+ / \lambda^D \hat{\alpha}^+) - \hat{p}^v, \quad r_0^- = \hat{p}^v - F^{-1}(\pi_0^- / \lambda^D \hat{\alpha}^-) \quad (6)$$

Taking downward BRS as an example, its demand characteristics is shown in Fig. 4. Given lower BRS prices, the producer will choose to purchase more. And the maximum amount of BRS it would like to buy is the headroom between DA schedule and the installed capacity.

The whole area enclosed by the demand curve and axes, namely, ($S_1+S_2+S_3$), represents the potential losses if no BRS is purchased. And, for another extreme case, if the price of BRS is 0, then the producers will purchase $(p^{v,max} - \hat{p}^v)$ downward BRS and can get rid of all the potential imbalance costs. For most of the situations, the producer will purchase BRS less than maximum amount and above zero. Therefore, the overall imbalance costs (OIC), which can be defined as the summation of imbalance penalties and payments for BRSs, can be represented by the area of ($S_2+S_3$), and the area of $S_1$ is the consumer surplus of the BRS owner. This is interesting because the expected OIC of a VG producer is tied to the price of capacity. With higher BRS prices, the OIC, i.e. the area of ($S_2+S_3$), will increase. In contrast, if the system is sufficient in capacity and the price of BRS is low, the VG producer can hedge its exposure to deviation penalties at low cost. This is reasonable and compliant to the rule of cost causation. In other words, the pay-off of VG producers is not only decided by their intrinsic level of uncertainty, but also related to the adequacy of re-dispatch capacity.

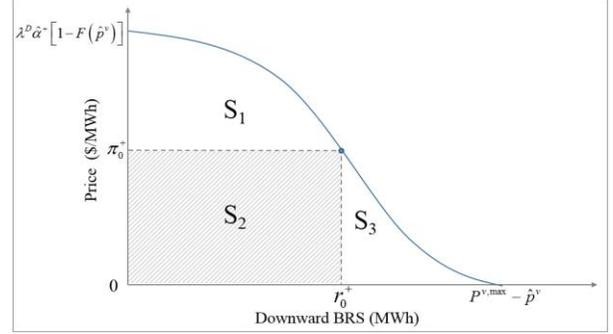

Fig. 4. Demand Characteristics of downward BRS of a VG producer.

### B. Supply

The characteristics of BRS supply is much more complex because it is related to many factors. Basically, the cost of providing BRS is the potential risk of deviating from the optimal scheduling point in the DA market. Without loss of generality, let us take a normal unit as an example. If there is no BRS contract, the revenue of the unit can be expressed as

$$R^g = \lambda^D \hat{p}^g + \left[ p^g (\lambda^R) - \hat{p}^g \right] \lambda^R \quad (7)$$

If the unit sells BRS to VG producers, and $\Delta r$ BRS is executed, the revenue of the unit will be changed to

$$R^{g*} = \lambda^D (\hat{p}^g + \Delta r) + \left[ p^g (\lambda^R) - \hat{p}^g - \Delta r \right] \lambda^R \quad (8)$$

where, positive and negative $\Delta r$ means executed upward and downward BRS respectively. After the DA market has cleared, $\lambda^D$ and $\hat{p}^g$ is given, while other quantities remain uncertain. Suppose that the owner of the unit is rational and bid its marginal cost to the market, then the execution of BRS will have no impact on the RT market clearing results. Therefore, the difference in pay-off is:

$$R^{g*} - R^g = (\lambda^D - \lambda^R) \Delta r \quad (9)$$

We assume that a single VG producer is a price taker and its deviation is independent from the difference between DA and RT market prices. Moreover, the expected $\Delta r$ and price gap are both close to zero. Therefore, the expectation of (9) is:

$$\mathbb{E}(R^{g*} - R^g) = \mathbb{E}(\lambda^D - \lambda^R) \mathbb{E}(\Delta r) \approx 0 \quad (10)$$

Hence, it seems that there is no extra cost of signing a BRS contract. However, equation (10) does represent an uncertain cash flow, which increases the risk of the BRS providers. More specifically, the increase in risk can be represented by the change in variance of cash flow:

$$\text{var}(R^{g*}) - \text{var}(R^g)$$
$$= \text{var}\left((\lambda^D - \lambda^R)\Delta r + [p^g(\lambda^R) - \hat{p}^g]\lambda^R\right) - \text{var}\left([p^g(\lambda^R) - \hat{p}^g]\lambda^R\right) \quad (11)$$

This quantity may be related to the type of units. For base-load units, the RT schedule is non-sensitive to RT prices since their marginal costs are usually below the clearing prices. Therefore, equation (12) can be approximately simplified as

$$\text{var}(R^{g*}) - \text{var}(R^g) = \text{var}\left((\lambda^D - \lambda^R)\Delta r\right) \quad (12)$$

However, for those marginal units or quasi-marginal units, this is not the case. Denote $\vartheta_1 = (\lambda^D - \lambda^R)\Delta r$ and $\vartheta_2 = [p^g(\lambda^R) - \hat{p}^g]\lambda^R$. $\vartheta_1$ and $\vartheta_2$ cannot be perfectly positively correlated. Combining (10) and (11), we have:

$$\text{var}(\vartheta_1 + \vartheta_2) - \text{var}(\vartheta_2) < \text{var}(\vartheta_1) \quad (13)$$

The left-hand side represents the incremental variance of marginal or quasi-marginal units, and the right-hand side corresponds to base-load units. Therefore, BRS contracts bring less incremental risk to marginal units than base-load units, which means that, with the same risk preferences, marginal units have advantages in providing BRS. Moreover, based-load units are more likely to be fully scheduled at DA market. Therefore, to offer upward BRS, they have to reduce its DA schedule to reserve enough capacity, which will incur extra opportunity costs. This is also possible for marginal units but is substantially less likely, meaning that marginal units face less opportunity costs to provide BRS. Additionally, it should be noted that if (10) does not hold, the expected difference can be estimated and added to the cost of offering BRS.

It might not be straightforward to derive an explicit BRS supply curve of a unit. But it is clear that the cost of BRS corresponds to the compensation of risk. With higher risk, higher return will be required. And, providers that are more sensitive to market prices, such as marginal units, storage facilities and elastic demands, are all potential suppliers to the BRS market.

## IV. IMPACTS OF NETWORK CONSTRAINTS ON BRS

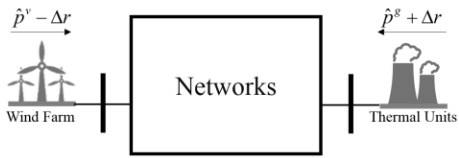

Fig. 5. Executed BRS contract may potentially change system power flow.

As shown in Fig. 5, a BRS contract enables the involved parties, i.e. a VG producer and a BRS provider, to switch DA scheduled output without being affected by market signals. If network congestion does exist, such privilege might hurt other participants in the system (cause infeasibility of DA power flow). In other words, market signals should be given to incent the VG producers to purchase re-dispatch capacity from "proper" locations.

To achieve that, the most convenient way is to forbid BRS transactions between congested regions. However, sometimes congestion can happen inside regions. According to the available methods for transmission fee and power losses allocation [13]-[14], a potential alternative is to carry out an adjustment economic dispatch process after BRS transactions are closed, guaranteeing that all the signed BRS contracts are physically feasible. Then, the incremental costs incurred by BRS contracts are allocated to BRS contract holders without affecting the DA settlements of other market participants. Due to limit of space, the detail method of allocating BRS congestion costs will be discussed in our future work.

## V. NUMERICAL CASES

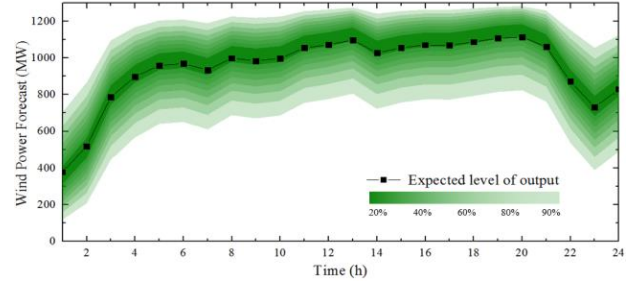

Fig. 6. DA probabilistic forecast of a specific wind farm.

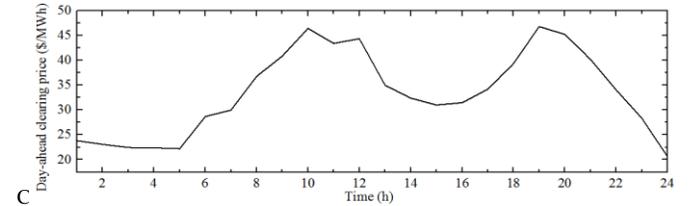

Fig. 7. DA energy clearing price.

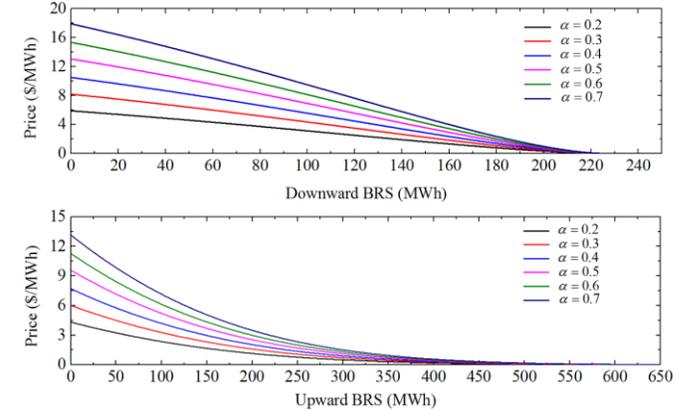

Fig. 8. Downward and upward BRS demand curve at hour 12 of the wind farm.

To demonstrate the demand characteristics of BRS, we take a hypothetical wind farm located in Texas as an example. The historical forecast and actual output of the wind farm was obtained from [15], based on which we derive the probabilistic forecast (see Fig. 6). It is assumed that the wind power follows a Beta distribution, and the forecast error variance is conditional to the predicted mean value. Fig. 7 demonstrates the DA market clearing price of a certain day obtained from ERCOT's website. For simplicity, suppose that $\alpha^+ = \alpha^- = \alpha$. Then, Fig. 8 presents the demand curve of downward and upward BRS at hour 12 with different penalty factors. The results indicate that, with larger penalty factors, heavier balancing responsibility is put on VG producers, and the demand for BRS at the same price level increases. In Fig. 9, the demand characteristics of BRS during the whole day with $\alpha$ set at 0.3 are demonstrated. A VG producer can seek to sign its optimal amount of BRS contracts according to the demand characteristics.

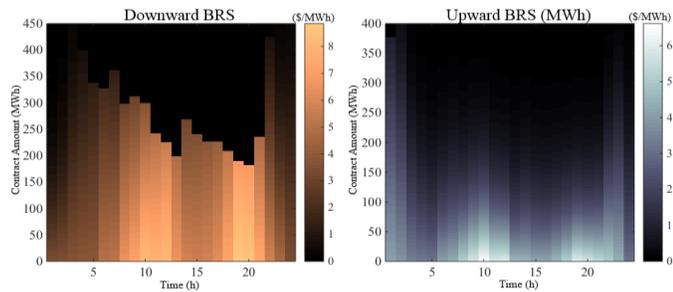

Fig. 9. BRS demand characteristics when penalty factors are set at 0.3.

To examine the impacts of BRS prices on VG producers, we simply assume that the BRS price is proportional to the DA energy price, and then adjust the ratio. Fig. 10 demonstrate the expected profits of the wind farm with different level of BRS prices, assuming that the producer has maximized its expected revenue by purchasing the optimal amount of BRS. The profit is defined as the expected revenue with BRS minus the costs of purchasing BRS contracts. Meanwhile, we modify the level of VG forecast uncertainty by changing the estimated variance of the distribution, and calculate the expected profit under different levels of uncertainty. The results are also plotted in Fig. 10. By comparison, we can find that a higher level of uncertainty, or worse forecast quality, can always lead to lower expected income. When the re-dispatch capacity is scarce and its price is high, the VG producers will tend to run with less or without BRS contracts. This is why the profits under each uncertainty level approach to a certain level as the price of BRS increases. When the system is sufficient in capacity and the price of BRS is low, the difference in expected profit is much less remarkable. This is reasonable since sufficient capacity makes the cost of uncertainty relatively low. In this way, the proposed BRS-based mechanism is effective in addressing both the level of uncertainty of VG producers and the system resource sufficiency.

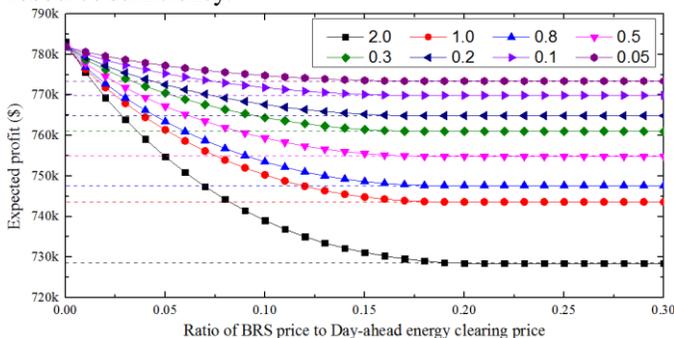

Fig. 10. Expected profits of a wind farm with different levels of BRS contracts for different levels of DA forecast uncertainty (variance). For example, the black curve symbolled by square is calculated by increasing the variance of forecast error to 2.0 times of the original value.

Due to the limited space, cases that demonstrate the entire market process of BRS and network congestion issues are not included.

## VI. CONCLUSIONS AND DISCUSSIONS

In this paper, we described a new product, BRS, and its market implementation to address the uncertainty of VG sources and the value of re-dispatch capacity. The basic idea is to allow the VG producers to purchase a certain amount of capacity after the DA scheduling process to avoid some of their deviation penalties. We analyzed the demand and supply characteristics of this product, and found that both sides have the incentives to participate in such transactions. Moreover, for VG producers, BRS can link their revenue to both their forecast quality and the availability of re-dispatch capacity resources. This is a desirable cost-causation feature that provides improves incentives to manage uncertainty in future power systems with large share of renewable energy.

For future work, the market mechanism considering impacts of network congestion on BRS transactions should be properly designed. More case studies on market clearing of such capacity products need to be done. Moreover, it should be carefully examined whether this product will lead to unexpected market manipulation behaviors.